\documentclass[11pt,draft]{article}
\usepackage{amsfonts}
\usepackage{amssymb}
\usepackage{amsmath}
\usepackage{amsthm}

\newtheorem{theorem}{Theorem}
\newtheorem{proposition}{Proposition}

\newtheorem{lemma}{Lemma}

\newtheorem{corollary}{Corollary}
\newcounter{xxx}

\renewcommand{\Re}{\mathop{\mathrm{Re}}}
\renewcommand{\Im}{\mathop{\mathrm{Im}}}

\newcommand{\supp}{\mathop{\mathrm{supp}}}
\newcommand{\card}{\mathop{\mathrm{card}}}
\newcommand{\clos}{\mathop{\mathrm{clos}}}
\newcommand{\const}{\mathop{\mathrm{const}}}
\newcommand{\GL}{\mathop{\mathrm{GL}}}

\newcommand{\SU}{\mathop{\mathrm{SU}}}
\newcommand{\SO}{\mathop{\mathrm{SO}}}

\newcommand{\cD}{\mathord{\mathcal{D}}}

\newcommand{\frC}{\mathord{\mathfrak{C}}}

\newcommand{\frh}{\mathord{\mathfrak{h}}}

\newcommand{\frU}{\mathord{\mathfrak{U}}}
\newcommand{\frV}{\mathord{\mathfrak{V}}}
\newcommand{\frg}{\mathord{\mathfrak{g}}}

\newcommand{\bbC}{\mathord{\mathbb{C}}}
\newcommand{\bbR}{{\mathord{\mathbb{R}}}}

\newcommand{\bbT}{\mathord{\mathbb{T}}}
\newcommand{\bbZ}{\mathord{\mathbb{Z}}}
\newcommand{\al}{{\mathord{\alpha}}}

\newcommand{\Ga}{{\mathord{\Gamma}}}
\newcommand{\vf}{{\mathord{\varphi}}}
\newcommand{\om}{{\mathord{\omega}}}
\newcommand{\si}{{\mathord{\sigma}}}

\newcommand{\la}{{\mathord{\lambda}}}
\newcommand{\La}{{\mathord{\Lambda}}}
\newcommand{\de}{{\mathord{\delta}}}
\newcommand{\De}{{\mathord{\Delta}}}

\newcommand{\ep}{{\mathord{\varepsilon}}}

\newcommand{\scal}[2]{\left<#1,#2\right>}

\let\ov=\overline
\let\td=\tilde

\let\dd=\partial

\date{}
\title{\bf A note on
common zeroes of Laplace--Beltrami eigenfunctions}
\author{V.M. Gichev}
\begin{document}
\maketitle
\begin{abstract}
Let $\De u+\la u=\De v+\la v=0$, where $\De$ is the
Laplace--Beltrami operator on a compact connected smooth manifold
$M$ and $\la>0$. If $H^1(M)=0$ then there exists $p\in M$ such
that $u(p)=v(p)=0$. For homogeneous $M$, $H^1(M)\neq0$ implies the
existence of a pair $u,v$ as above that has no common zero.
\end{abstract}
\section{Introduction}
Let $M$ be a compact connected closed orientable $C^\infty$-smooth
Riemannian $d$-dimensional manifold and
$\De$ be the Laplace--Beltrami operator on it. Set
\begin{eqnarray*}
E_\la=\{u\in C^2(M):\,\De u+\la u=0\}.
\end{eqnarray*}
The eigenspace $E_\la$ can be nontrivial only for $\la\geq0$.
If the contrary is not stated explicitly, we assume that
functions are real valued and linear spaces are finite dimensional;
$H^p(M)$ denotes de Rham cohomologies.
\begin{theorem}\label{comze}
Let $M$ be as above.
\begin{itemize}
\item[\rm (1)] Suppose $H^1(M)=0$. Then for any $\la\neq0$ and
each pair $u,v\in E_\la$ there exists $p\in M$ such that
$u(p)=v(p)=0$. \item[\rm (2)] If $M$ is a homogeneous space of a
compact Lie group of isometries then the converse is true:
$H^1(M)\neq0$ implies the existence of $\la\neq0$ and a pair
$u,v\in E_\la$ without common zeroes.
\end{itemize}
\end{theorem}
The circle $\bbT=\bbR/2\pi\bbZ$ and functions $u(t)=\cos t$,
$v(t)=\sin t$ provide the simplest example for (2).  Moreover, (2)
is an easy consequence of this example and the following
observation: for homogeneous Riemannian manifolds $M=G/H$, where
$G$ is compact and connected, $H^1(M)\neq0$ is equivalent to the
existence of $G$-equivariant mapping $M\to\bbT$ for some
nontrivial action of $G$ on $\bbT$.

The corollary below gives the answer to the question in \cite{GHV}:
is it true that each orbit of a compact connected irreducible linear group,
acting in a complex vector space, meets any hyperplane?
I am grateful to P. de la Harpe for making me aware of this question
which in fact was the starting point for this note.

\begin{corollary}\label{hyper}
Let $V$ be a complex linear space, $\dim V>1$, and $G\subset\GL(V)$
be a compact connected irreducible group. Then for any $v\in V$ and
every linear subspace $H\subset V$ of complex codimension $1$ there exists
$g\in G$ such that $gv\in H$.
\end{corollary}
There is a real version of this corollary.
Let $\tau$ be a real linear irreducible representation of a compact
connected Lie group $G$ in a real linear space $V_\tau$.
We may assume that $V_\tau$ is endowed with the invariant inner product
$\scal{\ }{\ }$ and that $G$ is equipped with a bi-invariant
Riemannian metric. Let $M_\tau$ be
the space of its matrix elements; by definition, $M_\tau$ is the linear span
of functions $t_{xy}(g)=\scal{\tau(g)x}{y}$, $x,y\in V_\tau$. Then either
$\tau$ admits an invariant complex structure or its complexification
is irreducible.
It follows from the Schur lemma that $\De u=\la_\tau u$ for each
$u\in M_\tau$, where $\De$ is the Laplace--Beltrami operator for
a bi-invariant metric on $G$. Let us fix $\De$ and denote by $\La_\si$,
where $\si$ is a finite dimensional real representation, the
spectrum of $\De$ on $M_\si$; it is the union
of $\la_\tau$ over all irreducible components $\tau$ of $\si$.
\begin{corollary}\label{codi2}
Let $G$ be a compact connected semisimple Lie group, $\si$ be
%its real linear finite dimensional representation
as above. Suppose that $\La_\si$ is a single point $\la\neq0$.
Then the orbit of any vector in $V_\si$ meets each linear subspace
of codimension $2$.\footnote{If $\La_\si$ is not a single point
then the assertion is not true. The simplest example is the
representation of $\SO(3)$ in the space of harmonic polynomials on
$\bbR^3$ of the type $l(x)+q(x)$, where $l$ is linear and $q$ is
quadratic. Let us fix $x_0\neq0$ and define a subspace $W$ of
codimension 2 by equalities $l(x_0)=0$ and $q(x_0)=0$. If $q_0$ is
nondegenerate then there exists $l_0$ such that $l^{-1}_0(0)\cap
q^{-1}_0(0)=\{0\}$; the orbit of $l_0+q_0$ does not intersect
$W$.}
\end{corollary}
If $u$ is an eigenfunction of $\De$ on a Riemannian manifold $M$ then
\begin{eqnarray*}
N_u=\{x\in M:\,u(x)=0\}
\end{eqnarray*}
is said to be the {\it nodal set}, and connected components of
its complement $M_u=M\setminus N_u$ are called {\it nodal domains}.
In the following lemma, we formulate  the main step in the proof of
the theorem {(the fact seems to be known but I failed to find a
reference)}.
\begin{lemma}\label{nesod}
Let $u,v\in E_\la$, $u,v\neq0$, and let $U$, $V$ be nodal domains for
$u$, $v$, respectively. If $U\subseteq V$ then $u=cv$ for some $c\in\bbR$.
\end{lemma}

There are many natural questions concerning the distribution
of common zeroes;
%.For example, is there an estimate from below for $(d-2)$-hausdorff
%measure for this set depending only on $\la$ and the geometry of $M$?
%These questions
they seem to be difficult. We prove a very particular result
for $d=\dim M=2$. Note that $M$ is diffeomorphic to the sphere $S^2$
if $d=2$ and $H^1(M)=0$.
\begin{proposition}\label{twopo}
Let $d=2$, $H^1(M)=0$, $\la\ne0$, $u\in E_\la$.
Suppose that zero is not a critical value for $u$.
Then for each $v\in E_\la$ every connected component of $N_u$
contains at least two points of $N_v$.
\end{proposition}
\noindent In fact, each component is a Jordan contour and supports
a positive measure which annihilates $E_\la$.

Let $M$ be the unit sphere $S^2\subset\bbR^3$ with the standard metric.
Then $\la_n=n(n+1)$ is $n$-th eigenvalue of $\De$. The corresponding
eigenspace $E_n=E_{\la_n}$ consists of spherical harmonics which can be
defined as restrictions to $S^2$ of harmonic (with respect to the ordinary
Laplacian in $\bbR^3$) homogeneous polynomials of degree $n$ in
$\bbR^3$; $\dim E_n=2n+1$. The space $E_n$ is spanned by zonal spherical
harmonics
$l_{a,n}(x)=L_n(\scal{x}{a})|_{S^2}$, where $a\in S^2$ and $L_n$ is the
Legendre polynomial. The nodal set for $l_{a,n}$ is the union of $n$ circles
$$\{x\in S^2:\,\scal{x}{a}=x_k\},$$
where $x_1,\dots,x_n\in[-1,1]$ are zeroes of $L_n$.
%(note that $l_{a,n}(x)$ is not harmonic on $\bbR^3$; it coincides with some
%harmonic polynomial on $S^2$).
Set $u=l_{a,n}$, $v=l_{b,n}$, $n(a,b)=\card(N_u\mathbin\cap N_v)$.
{Projections of $N_u$ and $N_v$ to the plane $\pi_{ab}$ passing
through $a$ and $b$ are families of segments in the unit disc in
$\pi_{ab}$ with endpoints in the unit circle which are orthogonal
to $a$ and $b$, respectively.
Outside the boundary circle, the preimage of each point is a
pair of points. Further, $N_u$, $N_v$ are symmetric with
respect to $\pi_{ab}$. Hence $N_u\mathbin\cap N_v$ corresponds to
the intersection of the segments. It makes possible to calculate or
estimate $n(a,b)$.
In particular, if $a$ and $b$ are sufficiently close then $n(a,b)=2n$;
if $a\perp b$ then $n(a,b)\approx cn^2$, where $c$ can be calculated
explicitly since zeroes of $L_n$ are distributed uniformly in $[-1,1]$.}
%\begin{eqnarray*}
%c=\frac1{\pi^2}
%\mathop{\int\int}\limits_{x^2+y^2<1}\frac{dx\,dy}{\sqrt{1-x^2}\,\sqrt{1-y^2}}.
%\end{eqnarray*}
The set $N_u\mathbin\cap N_v$ can be infinite for independent
$u,v\in E_n$, for instance, it can be a big circle or a family of
parallel circles in $S^2$ { (this is true for suitable $u,v$ of
the type $P(\cos\theta)\cos(k\vf+\al)$, where $P$ is a polynomial,
$\theta,\vf$ are Euler coordinates in $S^2$, $k=1,\dots,n$)}. I do
not know if there are other nontrivial examples of infinite sets
$N_u\mathbin\cap N_v$ as well as examples of $u,v\in E_n$ such
that $\card(N_u\mathbin\cap N_v)<2n$\footnote{with multiplicities,
or for generic $u,v$. If $n=2$ then $4\leq\card(N_u\mathbin\cap
N_v)\leq 8$; for $n=1$, $\card(N_u\mathbin\cap N_v)=2$.}.

{It is natural to ask if something like Theorem~\ref{comze} is
true for three or more eigenfunctions. Here is an example.} Let
$S^3$ be the unit sphere in $\bbC^2$ and set $u=|z_1|^2-|z_2|^2$,
$v=\Re z_1\ov z_2$, $w=\Im z_1\ov z_2$. These three
Laplace--Beltrami eigenfunctions have no common {zeroes} in $S^3$.
They are matrix elements of { the} three dimensional
representation of $\SU(2)\cong S^3$ and correspond to three linear
functions on $S^2\subset\bbR^3$; the homogeneous space $M$ admits
an equivariant mapping $M\to S^2$. { Perhaps, the latter property
could be the right replacement of the assumption $H^1(M)\neq0$ in
a version of Theorem~\ref{comze} for homogeneous spaces and three
eigenfunctions.}

\section{Proof of results}
By $\rho$ we denote the Riemannian metric in $M$, $\om$ is the
volume $n$-form. The metric $\rho$ identifies tangent and
cotangent bundles, hence it extends to $T^*M$.
%\footnote{and to its exterior powers.}
Let $D$ be a domain in $M$, $C^2_c(D)$ be the set of all functions
in $C^2(D)$ with compact support in $D$, $W_0$ be the closure of
$C^2_c(D)$ in the Sobolev class $W^1_2(D)$ which consist of
functions whose first derivatives (in the sense of the
distribution theory) are square integrable functions. There is the
natural unique up to equivalence norm making it a Banach space.
For all $u,v\in C^2_c(D)${\ }\footnote{thus, $\De=-(d\de+\de d)$,
where $\de$ is the adjoint operator for $d$; due to the choice of
the sign, $\De$ is the ordinary Laplacian in the Euclidean case.}
\begin{eqnarray*}
\int_D\rho(du,\,dv)\,\om=-\int_D u\De v\,\om=-\int_D v\De u\,\om.
\end{eqnarray*}
%(this equality could be considered as a definition of $\De$).
Hence for every $u,v,w\in C^2_c(D)$
\begin{eqnarray}\label{bpart}
\int_D u\rho(dv,\,dw)\,\om=-\int_D v(\rho(du,\,dw)+u\De w)\,\om.
\end{eqnarray}
For a domain $D\subseteq M$ and a function $u\in W^1_2(D)$, let
\begin{eqnarray*}
\cD_D(u)=\int_D\rho(du,\,du)\,\om
\end{eqnarray*}
be the Dirichlet form. {In most cases}, we shall omit the index.
For the sake of completeness, we give a proof of the classical result
which states that a positive eigenfunction corresponds to the first
eigenvalue which is multiplicity free. The proof follows
\setcounter{xxx}{6}\cite[Ch. \Roman{xxx}, \S7]{CH}.
\begin{lemma}\label{class}
Let $D$ be a domain in $M$, $v\in C^2(D)\mathbin\cap W_0$.
Suppose $v>0$ and $\De v+\la v=0$ on $D$. Then for all $u\in W_0$
\begin{eqnarray}\label{estif}
\cD(u)\geq\la\int_{D}u^2\,\om,
\end{eqnarray}
and the equality holds if and only if $u=cv$ in $D$ for some $c\in\bbR$.
\end{lemma}
\begin{proof}
Since $v>0$, each $u\in C^2_c(D)$ admits the unique factorization
$u=\eta v$, where $\eta\in C^2_c(D)$. Due to {\rm(\ref{bpart})}
and the equality $2\eta v\rho(d\eta,\,dv)=v\rho(d\eta^2,\,dv)$,
\begin{eqnarray*}
&\!\!\!\!\cD(u)=&\int_D\rho(d(\eta v),\,d(\eta v))\,\om=\\
&&\int_D\left(v^2\rho(d\eta,\,d\eta)+2\eta v\rho(d\eta,\,dv)+
\eta^2\rho(dv,\,dv)\right)\,\om=\\
&&\int_D\left(v^2\rho(d\eta,\,d\eta)+\eta^2\rho(dv,\,dv)\right)\,\om
-\int_D\eta^2\left(\rho(dv,\,dv)+v\De v\right)\,\om=\\
&&\int_D\left(v^2\rho(d\eta,\,d\eta)+\la\eta^2 v^2\right)\,\om\geq
\la\int_D\eta^2v^2\,\om=\la\int_D u^2\,\om.
\end{eqnarray*}
Using the approximation, we get {\rm(\ref{estif})}.
Suppose that the equality in {\rm(\ref{estif})} holds for some $u\in W_0$.
Let $\eta_n$ be such that $\eta_nv\to u$ in ${W_0}$ as $n\to\infty$.
Then $\cD(\eta_n v)\to\cD(u)$. Due to the calculation above,
\begin{eqnarray*}
\lim_{n\to\infty}\int_Dv^2\rho(d\eta_n,\,d\eta_n)\,\om=0.
\end{eqnarray*}
Let $D'\subset D$ be a domain whose closure is contained in $D$.
Standard arguments show that any limit point of the sequence
$\{\eta_n\}$ in $W_2^1(D')$ is a constant function.
%\footnote{since $d\eta=0$ for any limit function $\eta$}.
Hence $u=cv$ in $D$ for some $c\in\bbR$. The converse is obvious.
\end{proof}
\begin{proof}[Proof of Lemma~\ref{nesod}]
Let $D=V\supseteq U$; we may assume $v>0$ in $D$ and $u>0$ in $U$.
{Let $\td u$ be zero outside $U$ and coincide with $u$ in $U$.}
Clearly,\footnote{$\td u$ can be approximated in $W_0$ by
functions $u_n=\max\{\ep_n,u\}-\ep_n$ in $U$, $u_n=0$ outside $U$,
where $\ep_n>0$ are regular values for $u$ and $\ep_n\to0$ as
$n\to\infty$ (note that $u\in C^1(M)$).} $v,\td u\in W_0$.
Furthermore,
\begin{eqnarray*}
\cD_D(\td u)=\cD_U(u)=\la\int_Uu^2\,\om=\la\int_D\td u^2\,\om.
\end{eqnarray*}
By Lemma~\ref{class}, $u=cv$ in $D$. To conclude the proof, we refer
to Aronszajn's unique continuation theorem \cite{Ar} which implies
$u=cv$ on $M$.
\end{proof}
\begin{proof}[Proof of Theorem~\ref{comze}]
1) Let $\frU$ and $\frV$ be families of nodal domains for $u$ and $v$,
respectively. The assumption $\la\neq0$ and the orthogonality
relations imply $M_u,M_v\neq M$. Obviously, $u$ and $v$ have
no common zeroes if and only if $\frC=\frU\mathbin\cup\frV$
is a covering:
\begin{eqnarray}\label{cover}
M=\bigcup_{W\in\frC}W.
\end{eqnarray}
It is sufficient to prove, assuming {\rm(\ref{cover})},
that there exists a closed 1-form on $M$ which is not exact.
The covering $\frC$ has following properties:
\begin{itemize}
\item[\rm (A)] sets in $\frU$ are pairwise disjoint, and the same
is true for $\frV$; \item[\rm (B)] nor $U\subseteq V$ neither
$U\supseteq V$ for every $U\in\frU$, $V\in\frV$.
\end{itemize}
The first is obvious, the second is a consequence of Lemma~\ref{nesod}.
Also, Lemma~\ref{nesod} implies that
\begin{eqnarray}\label{intes}
U\mathbin\cap N_v\neq\emptyset,\quad V\mathbin\cap N_u\neq\emptyset\quad
\text{for all}\ U\in\frU,\, V\in\frV.
\end{eqnarray}
Due to {\rm(\ref{intes})}, $\frC$ is finite: $\frV$ covers the
compact set $N_u$ by open disjoint sets, and the same is true for
$\frU$ and $N_v$. This also means that a connected component $X$
of $N_u$ is contained in some nodal domain for $v$. Further, $u$
cannot keep its sign near $X$. Otherwise, we get a contradiction
assuming $u>0$ and applying the Green formula to  functions $u,1$
and the component of the set $u<\varepsilon$ which contains $X$
(for sufficiently small regular $\varepsilon>0$) \footnote{Indeed,
$\int_{\partial U_\ep}\de(u\,\om)=\int_{U_\ep}d(\de(
u\,\om))=-\int_{U_\ep}\De u\,\om=\la\int_{U_\ep}u\om>0$, where the
operator $\de$ is adjoint to $d$ and $U_\ep\supset X$ is the
component above. On the other hand, $\de(u\om)=-\frac{\partial
u}{\partial n}\om_\ep$, where $\frac{\partial }{\partial n}$ is
the outer normal and $\om_\ep$ is the volume form for $\partial
U_\ep$. Since $\frac{\partial u}{\partial n}\geq0$ on $\partial
U_\ep$, we get a contradiction.}
%Also, this property is local and can be deduced from the
%Euclidean case by the approximation.}.
Hence $X$ lies in the boundary of at least two domains in
$\mathfrak U$. Components of $N_v$ have this property with respect
to $\mathfrak V$. Let $\Gamma$ be the incidence graph for
$\mathfrak C$ whose family of vertices is $\mathfrak C$ and edges
join sets with nonempty intersection. For $\Gamma$, the conditions
above read as follows:
\begin{itemize}
\item [(a)] each edge {of $\Ga$} joins $\frU$ and $\frV$; \item
[(b)] any vertex is common for (at least) two different
edges\footnote{Otherwise, there exist $U\in\frU$, $V\in\frV$ such
that $U\subseteq\clos V$ or $V\subseteq\clos U$. If
$U\subseteq\clos V=V\cup\partial V$ then either $U\subseteq V$ or
$U\cap\partial V\neq\emptyset$. The first contradicts to (B), the
second implies the existence of a component $X$ of $N_v$ such that
$v$ keeps its sign near $X$.}.
\end{itemize}
It follows that $\Ga$ contains a nontrivial cycle $C$. Let $U\in\frU$
and $V\in\frV$ be consecutive vertices of $C$,
$Q=\partial U\mathbin\cap V$. Since $Q\mathbin\cap\dd V=\emptyset$
due to {\rm(\ref{cover})}, both $Q$ and
$\partial U\setminus Q=\partial U\setminus V$ are compact.
Hence there exists a smooth function $f$ on $M$ such that $f=1$
in a neighbourhood of $Q$ and $f=0$ near $\partial U\setminus Q$.
Then  $df=0$ on $\partial U$ and the 1-form
$\eta$ which is zero outside $U$ and coincides with $df$ on $U$
is {well} defined and smooth.  Obviously, $\eta$ is closed;
we claim that $\eta$ cannot be exact.

Suppose $\eta=dF$.  Then $F=\const$ on each
connected set which does not intersect $\supp\eta\subset U$.
Let $U_1=U,V_1=V,\dots,U_m,V_m$ be the cycle $C$. Then $m>1$
and we may assume that
\begin{eqnarray}\label{empin}
V_k\mathbin\cap U=\emptyset\quad\text{for}\ 1<k<m
\end{eqnarray}
replacing $C$ by a shorter cycle if necessary\footnote{We may
assume that $C$ contains no proper subcycle.}. If a curve in $V$
starts outside $U$ and comes into $U$ then it meets $U$ at a point
of $Q$. Hence there exists a curve $c_1$ in $V\setminus U$ with
endpoints in $Q$ and $U_2$. Analogously, there is a curve $c_2$
inside $V_m\setminus U$ which joins a point in $U_m$ with a point
in $V_m\mathbin\cap\partial U$. The set
\begin{eqnarray*}
X=c_1\mathbin\cup U_2\mathbin\cup V_2\mathbin\cup\dots\mathbin\cup
U_m\mathbin\cup c_2
\end{eqnarray*}
is connected; by A) and {\rm(\ref{empin})}, $X\mathbin\cap U=\emptyset$.
Therefore, $F$ is constant on $X$.  This contradicts to the choice of $f$
since the closure of $X$ has common points with $Q$ and
$\partial U\setminus Q$ (recall that $dF=df$ on $U$ and that $f$ takes
different values on these sets).

2) Let $M=G/H$, where $G$ is a compact group of isometries.
Since $M$ is connected, the identity component of $G$ acts
on $M$ transitively. Hence we may assume that $G$ is connected.
If $H^1(M)\neq0$ then there exists an invariant closed 1-form $\eta$
on $M$ that is not exact.
It can be lifted to the left invariant closed 1-form $\td\eta$ on the
universal covering group $\td G$.
Since $\td\eta$ is left invariant and closed (hence exact),
$\td\eta=d\chi$ for some nontrivial additive character $\chi:\,\td G\to\bbR$.
According to structure theorems, $\td G=\td S\times\bbR^k$, where $\td S$ is
compact, simply connected, and semisimple. Hence $\chi=0$ on $\td S$.
Further, $\eta$ is locally exact on $M$; thus $\chi=0$ on the preimage
$\td H$ of the group $H$ in $\td G$. Since $\td S$ is compact and normal,
$\td L=\td S\td H$ is a closed subgroup of $\td G$.
It follows that $\dim\td L<\dim G$. Let $S,L$ be subgroups
of $G$ which are covered by $\td S,\td L$, respectively. Thus, $L=SH$
is a closed proper subgroup of $G$.
The natural mapping $M=G/H\to G/L\cong\bbT^m$
can be continued to a nontrivial equivariant one $M\to\bbT$. Realizing
$\bbT$ as the unit circle in $\bbC$ we get a nonconstant function
whose real and imaginary parts satisfies the theorem.
\end{proof}
Conditions A) and B) imply $H^1(\frC,A)\neq0$ for any nontrivial
abelian group $A$. Indeed, every three {distinct} sets in $\frC$ {have}
empty intersection whence any function on the set of edges of $\Ga$
is a cocycle, {while} for each coboundary the sum of its values
along every cycle in $\Ga$
is zero. Besides, for homogeneous spaces of connected compact Lie groups
the condition $H^1(M)=0$ is equivalent to each of following ones:
$H_1(M,\bbZ)$ is finite;
$\pi_1(M)$ is finite; $\td M$ is compact; the semisimple part of $G$
acts on $M$ transitively; $\frg=\frh+[\frg,\frg]$, where $\frg$ and
$\frh$ are Lie algebras of $G$ and $H$, respectively. We omit the proof
which is easy.
\begin{proof}[Proof of Corollary~\ref{codi2}]
Since $G$ is semisimple, $H^1(G)=0$. Let $L\subset V_\si$ be a subspace
of codimension $2$, $x\in V_\si$, $y$ and $z$ be a linear base of $L^\bot$;
set $u(g)=\scal{\si(g)x}{y}$, $v(g)=\scal{\si(g)x}{z}$.
By Theorem~\ref{comze}, $u$ and $v$ have a common zero $g\in G$; then
$\si(g)x\perp y,z$ but this is equivalent to $\si(g)x\in L$.
\end{proof}
\begin{proof}[Proof of Corollary~\ref{hyper}]
Clearly, the semisimple part of $G$ is irreducible.
Hence $G$ can be assumed to be semisimple.
The condition $\dim V>1$ implies that the representation is not trivial.
Therefore, $\La_\tau$ is a single point $\la_\tau\neq0$ and the hyperplane
$H$ has the real codimension $2$ in $V$. Thus we may apply
Corollary~\ref{codi2}.
\end{proof}
{Note that the centre of $G$ consists of scalar matrices;
hence, if $G$ is not semisimple then
$H\mathbin\cap Gv$ includes $\bbT v$ for any $v\in H$, where
$\bbT$ is the unit circle in $\bbC$. Therefore, $H\mathbin\cap Gv$
is infinite for any $v\in V\setminus\{0\}$ in this case.}

In what follows, we assume that $M$ is diffeomorphic to
the sphere $S^2$. Let $D$ be a domain in $M$ bounded by a finite number
of smooth curves.
Then there exists a vector field $\frac{\dd }{\dd n}$ on $\dd D$
orthogonal to $\dd D$ such that the Green formula holds:
\begin{eqnarray}\label{greens}
\int_{\dd D}\left(u\frac{\dd v}{\dd n}-v\frac{\dd u}{\dd n}\right)\,ds=
\int_{D}\left(u\De v-v\De u\right)\,dm
\end{eqnarray}
for all smooth $u,v$, where $s$ and $m$ are linear and area measures
defined by $\rho$ on $\dd D$ and $D$, respectively. The vector field
$\frac{\dd }{\dd n}$ depends only on the local geometry of $\dd D$
and does not vanish on $\dd D$.
% (it can be easily derived from the usual Green formula).
\begin{lemma}\label{ortho}
Let $\la\ne0$, $u\in E_\la$, and  $C$ be a component of $N_u$.
If $C$ is a Jordan contour {that} contains no critical
points of $u$ then {there exists a strictly positive continuous
function $q$ on $C$ such that} $\int_{C}vq\,ds=0$ for all $v\in E_\la$.
\end{lemma}
\begin{proof}
Since $du\neq0$ on $C$, { it} is a smooth curve. Let $D$ be one of
the two domains bounded by $C$ due to Jordan Theorem.
Applying {\rm(\ref{greens})} to it, we get
\begin{eqnarray*}%\label{annih}
\int_{C}v\,\frac{\dd u}{\dd n}\,ds=0.
\end{eqnarray*}
{{Clearly}, if $\frac{\dd u}{\dd n}(p)=0$ for $p\in C$ then $p$
is a critical point. Hence either $q=\frac{\dd u}{\dd n}$
or $q=-\frac{\dd u}{\dd n}$  satisfies the lemma.}
\end{proof}
Proposition~\ref{twopo} is an easy consequence of Lemma~\ref{ortho}
(it remains to note that each component of $N_u$ is a Jordan contour
if $0$ is not a critical value).

\vskip1cm
V.M. Gichev\\
gichev@iitam.omsk.net.ru\\
Omsk Branch of \\
Sobolev Institute of Mathematics\\
Pevtsova, 13, 644099\\
Omsk, Russia

\begin{thebibliography}{9999}
\bibitem{Ar}
{Aronszajn N.\/}, {\it A unique continuation theorem for solutions of
elliptic partial differential equations of second order},
J. Math. Pures Appl., 36 (1957), 235-239.
\bibitem{CH}
{Courant R., Hilbert D.\/}, {\it Methoden der Mathematischen Physik},
Berlin, Verlag von Julius Springer, 1931.
\setcounter{xxx}{23}
\bibitem{GHV}
{Galindo J., de la Harpe P., Vust T.\/}, {\it Two Observations on
Irreducible Representations of Groups}, J. of Lie Theory, 12 (2002),
535--538.
\end{thebibliography}
\end{document}